\documentclass[11pt, twoside]{article}
\usepackage{latexsym}
\usepackage{amsmath}
\usepackage{amssymb}
\usepackage[all]{xy}
\usepackage{amsfonts}
\usepackage{verbatim}
\usepackage{amsthm}
\usepackage{mathrsfs}
\usepackage{epsfig}
\usepackage{xy}
\usepackage{array}
\usepackage{stmaryrd}
\usepackage{graphicx,color}
\usepackage{xcolor}
\usepackage{tikz}
\usetikzlibrary{arrows,calc}
\usepackage{etex}
\usepackage{mathdots}
\usepackage{float}
\usepackage{graphics}
\usepackage{pdflscape}
\usepackage{mathrsfs}

\usepackage{anysize,hyperref}
\input xypic
\xyoption{all}

\usepackage[perpage,symbol]{footmisc}
\usepackage{setspace}
\def\A{\mathcal{A}}

\def\P{\mathcal{P}}
\def\I{\mathcal{I}}

\def\C{\mathscr{C}}
\def\E{\mathbb{E}}
\def\del{\delta}
\def\s{\mathfrak{s}}
\def\id{\mathrm{id}}
\def\op{^\mathrm{op}}
\newcommand{\Xd}{\langle X_{\bullet},\del\rangle}  
\newcommand{\Yr}{\langle Y_{\bullet},\rho\rangle}  
\def\Ab{\mathit{Ab}}
\def\del{\delta}
\def\dr{\ar@{->}[r]}

\def\X{\mathscr{X}}

\newcommand{\CC}{{\bf{C}}^{n+2}_{\C}}\newcommand{\co}{\colon}
\newcommand{\uas}{^{\ast}}            
\newcommand{\sas}{_{\ast}}\newcommand{\ov}{\overset}
\newcommand{\lra}{\longrightarrow}
\newcommand{\ush}{^\sharp}           
\newcommand{\ssh}{_\sharp}

\usepackage{enumitem}
\begin{document}
\baselineskip=15pt
\title{\Large{\bf Higher Auslander-Reiten sequences revisited }}
\medskip
\author{Jian He, Hangyu Yin and Panyue Zhou\footnote{Corresponding author. Jian He was supported by the National Natural Science Foundation of China (Grant No. 12171230) and Youth Science and Technology Foundation of Gansu Provincial (Grant No. 23JRRA825). Panyue Zhou was supported by the National Natural Science Foundation of China (Grant No. 12371034) and the Hunan Provincial Natural Science Foundation of China (Grant No. 2023JJ30008).}}

\date{}

\maketitle
\def\blue{\color{blue}}
\def\red{\color{red}}

\newtheorem{theorem}{Theorem}[section]
\newtheorem{lemma}[theorem]{Lemma}
\newtheorem{corollary}[theorem]{Corollary}
\newtheorem{proposition}[theorem]{Proposition}
\newtheorem{conjecture}{Conjecture}
\theoremstyle{definition}
\newtheorem{definition}[theorem]{Definition}
\newtheorem{question}[theorem]{Question}
\newtheorem{remark}[theorem]{Remark}
\newtheorem{remark*}[]{Remark}
\newtheorem{example}[theorem]{Example}
\newtheorem{example*}[]{Example}
\newtheorem{condition}[theorem]{Condition}
\newtheorem{condition*}[]{Condition}
\newtheorem{construction}[theorem]{Construction}
\newtheorem{construction*}[]{Construction}

\newtheorem{assumption}[theorem]{Assumption}
\newtheorem{assumption*}[]{Assumption}

\baselineskip=17pt
\parindent=0.5cm

\begin{abstract}
\baselineskip=16pt
Let $(\C,\E,\s)$ be an $n$-exangulated category with enough projectives and enough injectives, and $\X$ be a cluster-tilting subcategory of $\C$.
Liu and Zhou have shown that the quotient category $\C/\X$ is an $n$-abelian category. In this paper, we prove that if $\C$ has Auslander-Reiten $n$-exangles, then $\C/\X$ has Auslander-Reiten $n$-exact sequences. Moreover, we also show that if a Frobenius $n$-exangulated category $\C$ has Auslander-Reiten $n$-exangles, then the stable category $\overline{\C}$ of $\C$ has Auslander-Reiten $(n+2)$-angles.\\[0.3cm]
\medskip
\end{abstract}

\pagestyle{myheadings}
\markboth{\rightline {\scriptsize  Jian He, Hangyu Yin and Panyue Zhou\hspace{2.5mm}}}
         {\leftline{\scriptsize   Higher Auslander--Reiten sequences revisited}}

\section{Introduction}
Auslander-Reiten theory, originally introduced by Auslander and Reiten in \cite{AR1,AR2}, has evolved into a foundational instrument for exploring the representation theory of Artin algebras. Subsequently, this theory has been generalized  to various contexts, including exact categories \cite{Ji}, triangulated categories \cite{H,RV} and their subcategories \cite{AS,J}, as well as specific additive categories \cite{L,J,S}, by various researchers.

The concept of extriangulated categories, introduced by Nakaoka and Palu in \cite{NP}, serves as a simultaneous generalization of both exact categories and triangulated categories. In this framework, a category is defined by a triplet $(\C, \E, \s)$, where $\C$ is an additive category, $\mathbb{E}: \C^{\rm op}\times \C \rightarrow {\rm Ab}$ is an additive bifunctor, and $\mathfrak{s}$ assigns to each $\delta\in \mathbb{E}(C,A)$ a class of $3$-term sequences with end terms $A$ and $C$, subject to certain axioms.
Recently, Herschend, Liu, and Nakaoka \cite{HLN} extended this concept to introduce the notion of $n$-exangulated categories for any positive integer $n$. This extension serves not only as a higher-dimensional analogue of extriangulated categories but also as a common generalization of $n$-exact categories, as defined by Jasso \cite{Ja}, and $(n+2)$-angulated categories, as defined by Geiss, Keller, and Oppermann \cite{GKO}. Interestingly, examples of $n$-exangulated categories exist that are neither $n$-exact nor $(n+2)$-angulated, as discussed in \cite{HLN, HLN1, LZ, HZZ2}.
Furthermore, the notion of Auslander-Reiten $n$-exact sequences in $n$-abelian categories was introduced by Fedele in \cite{F}. He, Hu, Zhang, and Zhou \cite{HHZZ} similarly delved into the concept of Auslander-Reiten $n$-exangles, exploring their existence in $n$-exangulated categories.

Cluster-tilting theory provides a systematic approach to constructing $n$-abelian categories from $n$-exangulated categories. To elaborate, Liu and Zhou \cite{LZ1} have shown that taking the quotient of an $n$-exangulated category by a cluster-tilting subcategory inherently results in an induced $n$-abelian structure on the quotient. This establishes a significant link between $n$-exangulated categories and $n$-abelian categories. In this paper, we delve into the relation between ``Auslander-Reiten sequences" within these two categories. Our first main result is the following.

\begin{theorem}{\rm (see Theorem \ref{main} for details)}
{\rm Let $(\C,\E,\s)$  be an $n$-exangulated category with split idempotents and let $\X$ be a cluster-tilting subcategory of $\C$.
Assume that $\C$ has enough projectives and enough injectives. If $\C$ has Auslander-Reiten $n$-exangles, then $\C/\X$ has Auslander-Reiten $n$-exact sequences.}
\end{theorem}

Liu and Zhou \cite{LZ} proved that the stable category of a Frobenius $n$-exangulated category is an $(n+2)$-angulated category, thereby establishing a significant connection between these categories. Building upon this foundational result, our second main result delves into exploring the relation between ``Auslander-Reiten sequences" within these two categories.

\begin{theorem}{\rm (see Theorem \ref{main1} for details)}
\rm Let $(\C,\E,\s)$  be a Frobenius $n$-exangulated category. If $\C$ has Auslander-Reiten $n$-exangles, then the stable category $\overline{\C}$ of $\C$ has Auslander-Reiten $(n+2)$-angles.
\end{theorem}

The structure of the paper is as follows: In Section 2, we offer a review of essential definitions and facts crucial to our study. In Section 3,
we provide proofs for the two main results of this paper.

\section{Preliminaries}

\subsection{$n$-exangulated categories}

let $\C$ be an additive category and $n$ be a positive integer. Suppose that $\C$ is equipped with an additive bifunctor $\E\colon\C\op\times\C\to{\rm Ab}$, where ${\rm Ab}$ is the category of abelian groups. Next we briefly recall some definitions and basic properties of $n$-exangulated categories from \cite{HLN}. We omit some
details here, but the reader can find them in \cite{HLN}.

{ For any pair of objects $A,C\in\C$, an element $\del\in\E(C,A)$ is called an {\it $\E$-extension} or simply an {\it extension}. We also write such $\del$ as ${}_A\del_C$ when we indicate $A$ and $C$. The zero element ${}_A0_C=0\in\E(C,A)$ is called the {\it split $\E$-extension}. For any pair of $\E$-extensions ${}_A\del_C$ and ${}_{A'}\del{'}_{C'}$, let $\delta\oplus \delta'\in\mathbb{E}(C\oplus C', A\oplus A')$ be the
element corresponding to $(\delta,0,0,{\delta}{'})$ through the natural isomorphism $\mathbb{E}(C\oplus C', A\oplus A')\simeq\mathbb{E}(C, A)\oplus\mathbb{E}(C, A')
\oplus\mathbb{E}(C', A)\oplus\mathbb{E}(C', A')$.

For any $a\in\C(A,A')$ and $c\in\C(C',C)$,  $\E(C,a)(\del)\in\E(C,A')\ \ \text{and}\ \ \E(c,A)(\del)\in\E(C',A)$ are simply denoted by $a_{\ast}\del$ and $c^{\ast}\del$, respectively.

Let ${}_A\del_C$ and ${}_{A'}\del{'}_{C'}$ be any pair of $\E$-extensions. A {\it morphism} $(a,c)\colon\del\to{\delta}{'}$ of extensions is a pair of morphisms $a\in\C(A,A')$ and $c\in\C(C,C')$ in $\C$, satisfying the equality
$a_{\ast}\del=c^{\ast}{\delta}{'}$.}
Then the functoriality of $\E$ implies $\E(c,a)=a_{\ast}(c^{\ast}\del)=c^{\ast}(a_{\ast}\del)$.
\begin{definition}\cite[Definition 2.7]{HLN}
Let $\bf{C}_{\C}$ be the category of complexes in $\C$. As its full subcategory, define $\CC$ to be the category of complexes in $\C$ whose components are zero in the degrees outside of $\{0,1,\ldots,n+1\}$. Namely, an object in $\CC$ is a complex $X_{\bullet}=\{X_i,d_i^X\}$ of the form
\[ X_0\xrightarrow{d_0^X}X_1\xrightarrow{d_1^X}\cdots\xrightarrow{d_{n-1}^X}X_n\xrightarrow{d_n^X}X_{n+1}. \]
We write a morphism $f_{\bullet}\co X_{\bullet}\to Y_{\bullet}$ simply $f_{\bullet}=(f_0,f_1,\ldots,f_{n+1})$, only indicating the terms of degrees $0,\ldots,n+1$.
\end{definition}
\begin{definition}\cite[Definition 2.11]{HLN}
By Yoneda lemma, any extension $\del\in\E(C,A)$ induces natural transformations
\[ \del\ssh\colon\C(-,C)\Rightarrow\E(-,A)\ \ \text{and}\ \ \del\ush\colon\C(A,-)\Rightarrow\E(C,-). \]
For any $X\in\C$, these $(\del\ssh)_X$ and $\del\ush_X$ are given as follows.
\begin{enumerate}
\item[\rm(1)] $(\del\ssh)_X\colon\C(X,C)\to\E(X,A)\ :\ f\mapsto f\uas\del$.
\item[\rm (2)] $\del\ush_X\colon\C(A,X)\to\E(C,X)\ :\ g\mapsto g\sas\delta$.
\end{enumerate}
We simply denote $(\del\ssh)_X(f)$ and $\del\ush_X(g)$ by $\del\ssh(f)$ and $\del\ush(g)$, respectively.
\end{definition}
\begin{definition}\cite[Definition 2.9]{HLN}
 Let $\C,\E,n$ be as before. Define a category $\AE:=\AE^{n+2}_{(\C,\E)}$ as follows.
\begin{enumerate}
\item[\rm(1)]  A pair $\Xd$ is an object of the category $\AE$ with $X_{\bullet}\in\CC$
and $\del\in\E(X_{n+1},X_0)$, called an $\E$-attached
complex of length $n+2$, if it satisfies
$$(d_0^X)_{\ast}\del=0~~\textrm{and}~~(d^X_n)^{\ast}\del=0.$$
We also denote it by
$$X_0\xrightarrow{d_0^X}X_1\xrightarrow{d_1^X}\cdots\xrightarrow{d_{n-2}^X}X_{n-1}
\xrightarrow{d_{n-1}^X}X_n\xrightarrow{d_n^X}X_{n+1}\overset{\delta}{\dashrightarrow}.$$
\item[\rm (2)]  For such pairs $\Xd$ and $\langle Y_{\bullet},\rho\rangle$,  $f_{\bullet}\colon\Xd\to\langle Y_{\bullet},\rho\rangle$ is
defined to be a morphism in $\AE$ if it satisfies $(f_0)_{\ast}\del=(f_{n+1})^{\ast}\rho$.

\end{enumerate}
\end{definition}
\begin{definition}\cite[Definition 2.13]{HLN}\label{def1}
 An {\it $n$-exangle} is an object $\Xd$ in $\AE$ that satisfies the listed conditions.
\begin{enumerate}
\item[\rm (1)] The following sequence of functors $\C\op\to\Ab$ is exact.
$$
\C(-,X_0)\xrightarrow{\C(-,\ d^X_0)}\cdots\xrightarrow{\C(-,\ d^X_n)}\C(-,X_{n+1})\xrightarrow{~\del\ssh~}\E(-,X_0)
$$
\item[\rm (2)] The following sequence of functors $\C\to\Ab$ is exact.
$$
\C(X_{n+1},-)\xrightarrow{\C(d^X_n,\ -)}\cdots\xrightarrow{\C(d^X_0,\ -)}\C(X_0,-)\xrightarrow{~\del\ush~}\E(X_{n+1},-)
$$
\end{enumerate}
In particular any $n$-exangle is an object in $\AE$.
A {\it morphism of $n$-exangles} simply means a morphism in $\AE$. Thus $n$-exangles form a full subcategory of $\AE$.
\end{definition}
\begin{definition}\cite[Definition 2.22]{HLN}
Let $\s$ be a correspondence which associates a homotopic equivalence class $\s(\del)=[{}_A{X_{\bullet}}_C]$ to each extension $\del={}_A\del_C$. Such $\s$ is called a {\it realization} of $\E$ if it satisfies the following condition for any $\s(\del)=[X_{\bullet}]$ and any $\s(\rho)=[Y_{\bullet}]$.
\begin{itemize}
\item[{\rm (R0)}] For any morphism of extensions $(a,c)\co\del\to\rho$, there exists a morphism $f_{\bullet}\in\CC(X_{\bullet},Y_{\bullet})$ of the form $f_{\bullet}=(a,f_1,\ldots,f_n,c)$. Such $f_{\bullet}$ is called a {\it lift} of $(a,c)$.
\end{itemize}
In such a case, we simple say that \lq\lq$X_{\bullet}$ realizes $\del$" whenever they satisfy $\s(\del)=[X_{\bullet}]$.

Moreover, a realization $\s$ of $\E$ is said to be {\it exact} if it satisfies the following conditions.
\begin{itemize}
\item[{\rm (R1)}] For any $\s(\del)=[X_{\bullet}]$, the pair $\Xd$ is an $n$-exangle.
\item[{\rm (R2)}] For any $A\in\C$, the zero element ${}_A0_0=0\in\E(0,A)$ satisfies
\[ \s({}_A0_0)=[A\ov{\id_A}{\lra}A\to0\to\cdots\to0\to0]. \]
Dually, $\s({}_00_A)=[0\to0\to\cdots\to0\to A\ov{\id_A}{\lra}A]$ holds for any $A\in\C$.
\end{itemize}
Note that the above condition {\rm (R1)} does not depend on representatives of the class $[X_{\bullet}]$.
\end{definition}
\begin{definition}\cite[Definition 2.23]{HLN}
Let $\s$ be an exact realization of $\E$.
\begin{enumerate}
\item[\rm (1)] An $n$-exangle $\Xd$ is called an $\s$-{\it distinguished} $n$-exangle if it satisfies $\s(\del)=[X_{\bullet}]$. We often simply say {\it distinguished $n$-exangle} when $\s$ is clear from the context.
\item[\rm (2)]  An object $X_{\bullet}\in\CC$ is called an {\it $\s$-conflation} or simply a {\it conflation} if it realizes some extension $\del\in\E(X_{n+1},X_0)$.
\item[\rm (3)]  A morphism $f$ in $\C$ is called an {\it $\s$-inflation} or simply an {\it inflation} if it admits some conflation $X_{\bullet}\in\CC$ satisfying $d_0^X=f$.
\item[\rm (4)]  A morphism $g$ in $\C$ is called an {\it $\s$-deflation} or simply a {\it deflation} if it admits some conflation $X_{\bullet}\in\CC$ satisfying $d_n^X=g$.
\end{enumerate}
\end{definition}

\begin{definition}\cite[Definition 2.32]{HLN}
An {\it $n$-exangulated category} is a triplet $(\C,\E,\s)$ of additive category $\C$, additive bifunctor $\E\co\C\op\times\C\to\Ab$, and its exact realization $\s$, satisfying the following conditions.

(EA1) Let $A\ov{f}{\lra}B\ov{g}{\lra}C$ be any sequence of morphisms in $\C$. If both $f$ and $g$ are inflations, then so is $g\circ f$. Dually, if $f$ and $g$ are deflations, then so is $g\circ f$.

(EA2) For $\rho\in\E(D,A)$ and $c\in\C(C,D)$, let ${}_A\langle X_{\bullet},c\uas\rho\rangle_C$ and ${}_A\Yr_D$ be distinguished $n$-exangles. Then $(\id_A,c)$ has a {\it good lift} $f_{\bullet}$, in the sense that its mapping cone gives a distinguished $n$-exangle $\langle M^f_{\bullet},(d^X_0)\sas\rho\rangle$.

(EA2$\op$) Dual of {\rm (EA2)}.

Note that the case $n=1$, a triplet $(\C,\E,\s)$ is a  $1$-exangulated category if and only if it is an extriangulated category, see \cite[Proposition 4.3]{HLN}.
\end{definition}
We recall the notion of Auslander-Reiten $n$-exangles in $n$-exangulated categories  from \cite{HHZZ}.

Let $\C$ be an additive category. We denote by ${\rm rad}_{\C}$ the Jacobson radical of $\C$, which is defined as follows:
$${\rm rad}_{\C}(A,B)=\{f\colon A\to B~|~~1_A-gf~~ \mbox{is invertible for any}~ g\colon B\to A\},$$
where $A,B\in\C$.
\begin{definition}\cite[Definition 3.1]{HLN} Let $\C$ be an $n$-exangulated category.
A distinguished $n$-exangle
$$A_0\xrightarrow{\alpha_0}A_1\xrightarrow{\alpha_1}A_2\xrightarrow{\alpha_2}\cdots\xrightarrow{\alpha_{n-2}}A_{n-1}
\xrightarrow{\alpha_{n-1}}A_n\xrightarrow{\alpha_n}A_{n+1}\overset{\delta}{\dashrightarrow}$$
in $\C$ is called an \emph{Auslander-Reiten $n$-exangle }if
$\alpha_0$ is left almost split, $\alpha_n$ is right almost split and
when $n\geq 2$, $\alpha_1,\alpha_2,\cdots,\alpha_{n-1}$ are in ${\rm rad}_{\C}$.
\end{definition}

\begin{definition}
Let $\C$ be an $n$-exangulated category. If for any  non-projective indecomposable object $A\in\C$, there exists an Auslander-Reiten $n$-exangle ending at $A$, and for any  non-injective indecomposable object $B\in\C$, there exists an Auslander-Reiten $n$-exangle starting at $B$.
In this case, we say that $\C$ has Auslander-Reiten $n$-exangles.
\end{definition}

\subsection{$n$-abelian categories}
Let $\A$ be an additive category and $f\colon A\rightarrow B$ a morphism in $\A$. A \emph{weak cokernel} of $f$ is a morphism
$g\colon B\rightarrow C$ such that for any $X\in\A$ the sequence of abelian groups
$$\A(C,X)\xrightarrow{~\A(g,X)~}\A(B,X)\xrightarrow{~\A(X,f)~}\A(A,X)$$
is exact. Equivalently, $g$ is a weak cokernel of $f$ if $gf=0$ and for each morphism
$h\colon B\rightarrow X$ such that $hf=0$ there exists a (not necessarily unique) morphism
$p\colon C\rightarrow X$ such that $h=pg$. Clearly, a weak cokernel $g$ of $f$ is a cokernel of $f$ if and only if $g$ is an epimorphism.
The concept of a \emph{weak kernel} is defined dually.

\begin{definition}\cite[Definiton 2.2]{Ja}
Let $\A$ be an additive category and $f_0\colon A_0\rightarrow A_1$  a morphism in
$\A$. An $n$-\emph{cokernel} of $f_0$ is a sequence
$$(f_1,f_2,\cdots,f_{n})\colon A_1\xrightarrow{~f_1~}A_2\xrightarrow{~f_2~}\cdots \xrightarrow
{~f_{n-1}~}A_n\xrightarrow
{~f_{n}~}A_{n+1}$$
such that the induced sequence of abelian groups
$$\xymatrix{0\xrightarrow{~~}\A(A_{n+1},B)\xrightarrow{~~} \A(A_{n},B)\xrightarrow{~~}\cdots \xrightarrow{~~}
\A(A_{1},B)\xrightarrow{~~}
\A(A_{0},B)}$$
is exact for each object $B\in\A$. That is, the morphism $f_i$ is a weak cokernel of $f_{i-1}$ for all $i=1,2,\cdots,n-1$ and $f_{n}$ is a cokernel of $f_{n-1}$. In this case, we say the sequence
$$A_0\xrightarrow{~f_0~}A_1\xrightarrow{~f_1~}A_2\xrightarrow{~f_2~}\cdots \xrightarrow
{~f_{n-1}~}A_n\xrightarrow
{~f_{n}~}A_{n+1} $$
is \emph{right $n$-exact}.

We can define \emph{$n$-kernel} and \emph{left $n$-exact} sequence dually. The sequence
$$A_0\xrightarrow{~f_0~}A_1\xrightarrow{~f_1~}A_2\xrightarrow{~f_2~}A_3\xrightarrow{~f_3~}\cdots\xrightarrow
{~f_{n-1}~}A_n\xrightarrow{~f_n~}A_{n+1}$$
is called \emph{$n$-exact} if it is
both right $n$-exact and left $n$-exact.
\end{definition}

\begin{definition}\cite[Definiton 3.1]{Ja}\label{def0}
Let $n$ be a positive integer. An \emph{$n$-abelian category} is an additive category $\A$
which satisfies the following axioms:
\begin{itemize}
\item[(A0)] The category $\A$ has split idempotents.

\item[(A1)] Every morphism in $\A$ has an $n$-kernel and an $n$-cokernel.

\item[(A2)] For every monomorphism $f_0\colon A_0\to A_1$ in $\A$ there exists an $n$-exact sequence:
$$A_0\xrightarrow{~f_0}A_1\xrightarrow{~f_1~}A_2
\xrightarrow{~f_2~}\cdots\xrightarrow{~f_3~}
A_{n-1}\xrightarrow{~f_{n-1}~}A_{n}
\xrightarrow{~f_n~}A_{n+1}.$$

\item[(A2)$^{\textrm{op}}$] For every epimorphism $g_n\colon B_n\to B_{n+1}$ in $\A$ there exists an $n$-exact sequence:
$$B_0\xrightarrow{~g_0}B_1\xrightarrow{~g_1~}B_2
\xrightarrow{~g_2~}\cdots\xrightarrow{~g_3~}
B_{n-1}\xrightarrow{~g_{n-1}~}B_{n}
\xrightarrow{~g_n~}B_{n+1}.$$
\end{itemize}
\end{definition}
We recall the notion of Auslander-Reiten $n$-exact sequences in $n$-abelian categories  from \cite{F}.
\begin{definition}\cite[Definition 4.3]{F} Let $\A$ be an $n$-abelian category.
An $n$-exact sequence
$$0\xrightarrow{} A_0\xrightarrow{\alpha_0}A_1\xrightarrow{\alpha_1}A_2\xrightarrow{\alpha_2}\cdots\xrightarrow{\alpha_{n-2}}A_{n-1}
\xrightarrow{\alpha_{n-1}}A_n\xrightarrow{\alpha_n}A_{n+1}\xrightarrow{}0$$
in $\A$ is called an \emph{Auslander-Reiten $n$-exact sequence }if
$\alpha_0$ is left almost split, $\alpha_n$ is right almost split and
when $n\geq 2$, $\alpha_1,\alpha_2,\cdots,\alpha_{n-1}$ are in ${\rm rad}_{\A}$.
\end{definition}

\begin{definition}
Let $\A$ be an $n$-abelian category. If for any  non-projective indecomposable object $A\in\A$, there exists an Auslander-Reiten $n$-exact sequence ending at $A$, and for any  non-injective indecomposable object $B\in\A$, there exists an Auslander-Reiten $n$-exact sequence starting at $B$.
In this case, we say that $\A$ has Auslander-Reiten $n$-exact sequences.
\end{definition}

\subsection{$(n+2)$-angulated categories}
Let $\C$ be an additive category and $\Sigma\colon \C\to\C$ be an automorphism
of $\C$.
An $(n+2)$-$\Sigma$-sequence in $\C$ is a sequence of objects and morphisms
$$A_0\xrightarrow{f_0}A_1\xrightarrow{f_1}A_2\xrightarrow{f_2}\cdots\xrightarrow{f_{n-1}}A_n\xrightarrow{f_n}A_{n+1}\xrightarrow{f_{n+1}}\Sigma A_0.$$
Its {\em left rotation} is the $(n+2)$-$\Sigma$-sequence
$$A_1\xrightarrow{f_1}A_2\xrightarrow{f_2}A_3\xrightarrow{f_3}\cdots\xrightarrow{f_{n}}A_{n+1}\xrightarrow{f_{n+1}}\Sigma A_0\xrightarrow{(-1)^{n}\Sigma f_0}\Sigma A_1.$$
A \emph{morphism} of $(n+2)$-$\Sigma$-sequences is  a sequence of morphisms $\varphi=(\varphi_0,\varphi_1,\cdots,\varphi_{n+1})$ such that the following diagram commutes
$$\xymatrix{
A_0 \ar[r]^{f_0}\ar[d]^{\varphi_0} & A_1 \ar[r]^{f_1}\ar[d]^{\varphi_1} & A_2 \ar[r]^{f_2}\ar[d]^{\varphi_2} & \cdots \ar[r]^{f_{n}}& A_{n+1} \ar[r]^{f_{n+1}}\ar[d]^{\varphi_{n+1}} & \Sigma A_0 \ar[d]^{\Sigma \varphi_0}\\
B_0 \ar[r]^{g_0} & B_1 \ar[r]^{g_1} & B_2 \ar[r]^{g_2} & \cdots \ar[r]^{g_{n}}& B_{n+1} \ar[r]^{g_{n+1}}& \Sigma B_0
}$$
where each row is an $(n+2)$-$\Sigma$-sequence. It is an {\em isomorphism} if $\varphi_0, \varphi_1, \varphi_2, \cdots, \varphi_{n+1}$ are all isomorphisms in $\C$.

\begin{definition}\cite[Definition 2.1]{GKO}
An $(n+2)$-\emph{angulated category} is a triple $(\C, \Sigma, \Theta)$, where $\C$ is an additive category, $\Sigma$ is an automorphism of $\C$ ($\Sigma$ is called $n$-suspension functor), and $\Theta$ is a class of $(n+2)$-$\Sigma$-sequences (whose elements are called $(n+2)$-angles), which satisfies the following axioms:
\begin{itemize}[leftmargin=3em]
\item[\bf (N1)] (a) The class $\Theta$ is closed under isomorphisms, direct sums and direct summands.

(b) For each object $A\in\C$ the trivial sequence
$$ A\xrightarrow{1_A}A\rightarrow 0\rightarrow0\rightarrow\cdots\rightarrow 0\rightarrow \Sigma A$$
belongs to $\Theta$.

(c) Each morphism $f_0\colon A_0\rightarrow A_1$ in $\C$ can be extended to $(n+2)$-$\Sigma$-sequence: $$A_0\xrightarrow{f_0}A_1\xrightarrow{f_1}A_2\xrightarrow{f_2}\cdots\xrightarrow{f_{n-1}}A_n\xrightarrow{f_n}A_{n+1}\xrightarrow{f_{n+1}}\Sigma A_0.$$

\item[\bf (N2)] An $(n+2)$-$\Sigma$-sequence belongs to $\Theta$ if and only if its left rotation belongs to $\Theta$.

\item[\bf (N3)] Each solid commutative diagram
$$\xymatrix{
A_0 \ar[r]^{f_0}\ar[d]^{\varphi_0} & A_1 \ar[r]^{f_1}\ar[d]^{\varphi_1} & A_2 \ar[r]^{f_2}\ar@{-->}[d]^{\varphi_2} & \cdots \ar[r]^{f_{n}}& A_{n+1} \ar[r]^{f_{n+1}}\ar@{-->}[d]^{\varphi_{n+1}} & \Sigma A_0 \ar[d]^{\Sigma\varphi_0}\\
B_0 \ar[r]^{g_0} & B_1 \ar[r]^{g_1} & B_2 \ar[r]^{g_2} & \cdots \ar[r]^{g_{n}}& B_{n+1} \ar[r]^{g_{n+1}}& \Sigma B_0
}$$ with rows in $\Theta$, the dotted morphisms exist and give a morphism of  $(n+2)$-$\Sigma$-sequences.

\item[\bf (N4)] In the situation of (N3), the morphisms $\varphi_2,\varphi_3,\cdots,\varphi_{n+1}$ can be chosen such that the mapping cone
$$A_1\oplus B_0\xrightarrow{\left(\begin{smallmatrix}
                                        -f_1&0\\
                                        \varphi_1&g_0
                                       \end{smallmatrix}
                                     \right)}
A_2\oplus B_1\xrightarrow{\left(\begin{smallmatrix}
                                        -f_2&0\\
                                        \varphi_2&g_1
                                       \end{smallmatrix}
                                     \right)}\cdots\xrightarrow{\left(\begin{smallmatrix}
                                        -f_{n+1}&0\\
                                        \varphi_{n+1}&g_n
                                       \end{smallmatrix}
                                     \right)} \Sigma A_0\oplus B_{n+1}\xrightarrow{\left(\begin{smallmatrix}
                                        -\Sigma f_0&0\\
                                        \Sigma\varphi_1&g_{n+1}
                                       \end{smallmatrix}
                                     \right)}\Sigma A_1\oplus\Sigma B_0$$
belongs to $\Theta$.
\end{itemize}
\end{definition}

We recall the notion of Auslander-Reiten $(n+2)$-angles in $(n+2)$-angulated categories  from \cite{F}.

\begin{definition}\cite[Definition 5.1]{F1}, \cite[Definition 3.8]{IY} Let $(\C, \Sigma, \Theta)$ be an $(n+2)$-angulated category.
An $(n+2)$-angle
$$A_0\xrightarrow{\alpha_0}A_1\xrightarrow{\alpha_1}A_2\xrightarrow{\alpha_2}\cdots\xrightarrow{\alpha_{n-2}}A_{n-1}
\xrightarrow{\alpha_{n-1}}A_n\xrightarrow{\alpha_n}A_{n+1}\xrightarrow{\alpha_{n+1}}\Sigma A_0$$
in $\C$ is called an \emph{Auslander-Reiten $(n+2)$-angle} if
$\alpha_0$ is left almost split, $\alpha_n$ is right almost split and
when $n\geq 2$, $\alpha_1,\alpha_2,\cdots,\alpha_{n-1}$ are in ${\rm rad}_{\C}$.
\end{definition}

\begin{definition}\cite[Theorem 3.8]{Z}
Let $(\C, \Sigma, \Theta)$ be an $(n+2)$-angulated category. If for any indecomposable object $C\in\C$, there exists an Auslander-Reiten $(n+2)$-angle ending at $C$, and for any indecomposable object $B\in\A$, there exists an Auslander-Reiten $(n+2)$-angle starting at $B$.
In this case, we say that $\C$ has Auslander-Reiten $(n+2)$-angles.
\end{definition}

\section{Main result}
We first recall the notion of \emph{cluster-tilting} subcategory from \cite{LZ1}.
\begin{definition}\label{defn}{\rm\cite[Definition 3.4]{LZ1}}
Let $\C$ be an $n$-exangulated category and $\X$ be a subcategory of $\C$.
$\X$ is called \emph{cluster-tilting} if

(1) $\E(\X,\X)=0$.

(2)  For any object $C\in\C$, there are two  distinguished $n$-exangles
$$X_0\xrightarrow{~~}X_1\xrightarrow{~~}\cdots\xrightarrow{~~}X_{n-1}\xrightarrow{~~}X_n\xrightarrow{~~}C\dashrightarrow$$
where $X_0, X_1,\cdots,X_{n}\in\X$ and
$$C\xrightarrow{~~}X'_1\xrightarrow{~~}X'_2\xrightarrow{~~}\cdots\xrightarrow{~~}X'_n\xrightarrow{~~}X'_{n+1}\dashrightarrow$$
where $X'_1,X'_2,\cdots,X'_{n+1}\in\X$.\end{definition}

Now we give some examples of cluster-tilting subcategories.

\begin{example}\label{ex2}
{\upshape Let $\Lambda$ be the algebra given by the following (infinity) quiver with relations $x^2=0$:
	 \begin{align}
	 	\begin{minipage}{0.6\hsize}
	 		\ \ \ \ \  \xymatrix{\begin{smallmatrix}1\end{smallmatrix}&\begin{smallmatrix}2\end{smallmatrix}\ar[l]_{x}
	 			&\begin{smallmatrix}3\end{smallmatrix}\ar[l]_{x}&\begin{smallmatrix}4\end{smallmatrix}\ar[l]_{x}&\begin{smallmatrix}\cdots\end{smallmatrix}\ar[l]_{x}&\begin{smallmatrix}n\end{smallmatrix}\ar[l]_{x}&\begin{smallmatrix}\cdots\end{smallmatrix}\ar[l]_{x}&\begin{smallmatrix}\end{smallmatrix}}\notag
	 	\end{minipage}
	 \end{align}
The Auslander-Reiten quiver of ${\rm mod}\Lambda$ is the following:}
$$\xymatrix@C=0.2cm@R0.2cm{
&&&\bullet \ar[dr] &&\bullet \ar[dr]&&\bullet \ar[dr] &&\bullet \ar[dr] &&\bullet \ar[dr] &&\bullet \ar[dr] &&\bullet \ar[dr] &&\bullet \ar[dr]&&\bullet\ar[dr]\\
&&\bullet\ar[ur]  &&\circ\ar[ur]  &&\spadesuit\ar[ur]  &&\circ \ar[ur]  &&\bullet\ar[ur] &&\circ\ar[ur] &&\spadesuit \ar[ur]  &&\circ\ar[ur]   &&\bullet\ar[ur]   &&\cdots
}
$$
where the object denoted by $\spadesuit$ and $\bullet $ appear periodically. Let $\C$ be the additive closure of all the indecomposable objects denoted by $\spadesuit$ and $\bullet $.
 Then $\C$ is a cluster-tilting subcategory of ${\rm mod}\Lambda$,
hence it is $2$-abelian (see \cite[Theorem 3.16]{Ja}). Let $\X$ be the additive closure of all the indecomposable objects denoted by $\bullet$.
It is straightforward to verify that  $\X$ is a cluster-tilting subcategory of $\C$.
\end{example}

\begin{example}
Consider $\Lambda$ as an $n$-representation finite algebra, and let $\mathscr{O}_{\Lambda}$ denote the $(n+2)$-angulated cluster category associated with $\Lambda$. Specifically, setting $n=3$, we define $\mathscr{T}=\mathscr{O}_{A^3_2}$, which represents the $5$-angulated (higher) cluster category of type $A_2$ (refer to \cite[Definition 5.2, Section 6, Section 8]{OT}). The indecomposable objects in this category correspond to the set
 $$\rm ind\mathscr{T}=\{1357, 1358, 1368, 1468, 2468, 2469, 2479, 2579, 3579 \}.$$
 The Auslander-Reiten quiver of $\mathscr{T}$ is depicted as follows:
$$
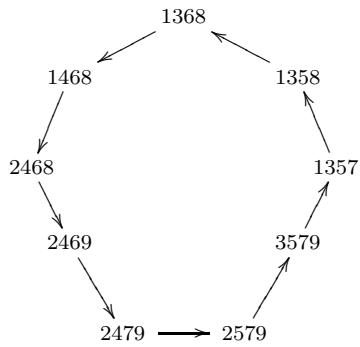
\begin{xy}
 (0,0)*+{\begin{smallmatrix}1368\end{smallmatrix}}="1",
(-15,-8)*+{\begin{smallmatrix}1468\end{smallmatrix}}="2",
(15,-8)*+{\begin{smallmatrix}1358\end{smallmatrix}}="3",
(-20,-20)*+{\begin{smallmatrix}2468\end{smallmatrix}}="4",
(20,-20)*+{\begin{smallmatrix}1357\end{smallmatrix}}="5",
(-15,-30)*+{\begin{smallmatrix} 2469\\
\end{smallmatrix}}="6",
(15,-30)*+{\begin{smallmatrix} 3579\\
\end{smallmatrix}}="7",
(-8,-42)*+{\begin{smallmatrix} 2479\\
\end{smallmatrix}}="8",
(8,-42)*+{\begin{smallmatrix} 2579\\
\end{smallmatrix}}="9",
\ar"1";"2", \ar"3";"1", \ar"2";"4", \ar"5";"3", \ar"4";"6",
\ar"7";"5",\ar"6";"8",\ar"9";"7",\ar"8";"9",
\end{xy}$$
It is easy to confirm that the subcategory
$${\rm add}T={\rm add}(1357\oplus1358\oplus1368\oplus1468)$$ is cluster-tilting.
\end{example}

Let $\C$ be an additive category and $\X$ a subcategory of $\C$.
Recall that we say a morphism $f\colon A \to B$ in $\C$ is an $\X$-\emph{monic} if
$$\C(f,X)\colon \C(B,X) \to \C(A,X)$$
is an epimorphism for all $X\in\X$. We say that $f$ is an $\X$-\emph{epic} if
$$\C(X,f)\colon \C(X,A) \to \C(X,B)$$
is an epimorphism for all $X\in\X$.
Similarly,
we say that $f$ is a left $\X$-approximation of $A$ if $f$ is an $\X$-monic and $ B\in\X$.
We say that $f$ is a right $\X$-approximation of $B$ if $f$ is an $\X$-epic and $A\in\X$.

\begin{definition} (\cite[Definition 3.1]{LZ})\label{dd1}
Let $(\C,\E,\s)$ be an $n$-exangulated category. A subcategory $\X$ of $\C$ is called
\emph{strongly contravariantly finite}, if for any object $C\in\C$, there exists a distinguished $n$-exangle
$$B\xrightarrow{}X_1\xrightarrow{}X_2\xrightarrow{}\cdots\xrightarrow{}X_{n-1}\xrightarrow{}X_{n}\xrightarrow{~g~}C\overset{}{\dashrightarrow}$$
where $g$ is a right $\X$-approximation of $C$ and $X_1,X_2,\cdots,X_n\in\X$.
Dually we can define \emph{strongly covariantly finite} subcategory. A strongly contravariantly finite and strongly  covariantly finite subcategory is called \emph{ strongly functorially finite}.
\end{definition}

\begin{remark}
In Definition \ref{defn}, it is clear that a cluster-tilting subcategory $\X$ is strongly functorially finite.
\end{remark}

\begin{definition}\label{def2}\cite[Definition 3.14 and Definition 3.15]{ZW}
Let $(\C,\E,\s)$ be an $n$-exangulated category.
\begin{itemize}
\item[(1)] An object $P\in\C$ is called \emph{projective} if, for any distinguished $n$-exangle
$$A_0\xrightarrow{\alpha_0}A_1\xrightarrow{\alpha_1}A_2\xrightarrow{\alpha_2}\cdots\xrightarrow{\alpha_{n-2}}A_{n-1}
\xrightarrow{\alpha_{n-1}}A_n\xrightarrow{\alpha_n}A_{n+1}\overset{\delta}{\dashrightarrow}$$
and any morphism $c$ in $\C(P,A_{n+1})$, there exists a morphism $b\in\C(P,A_n)$ satisfying $\alpha_n\circ b=c$.
We denote the full subcategory of projective objects in $\C$ by $\P$.
Dually, the full subcategory of injective objects in $\C$ is denoted by $\I$.

\item[(2)] We say that $\C$ {\it has enough  projectives} if
for any object $C\in\C$, there exists a distinguished $n$-exangle
$$B\xrightarrow{\alpha_0}P_1\xrightarrow{\alpha_1}P_2\xrightarrow{\alpha_2}\cdots\xrightarrow{\alpha_{n-2}}P_{n-1}
\xrightarrow{\alpha_{n-1}}P_n\xrightarrow{\alpha_n}C\overset{\delta}{\dashrightarrow}$$
satisfying $P_1,P_2,\cdots,P_n\in\P$. We can define the notion of having \emph{enough injectives} dually.
\end{itemize}
\end{definition}

Our main result is the following, which establish a connection between  Auslander-Reiten $n$-exangles in $n$-exangulated categories and Auslander-Reiten $n$-exact sequences in $n$-abelian categories.

\begin{theorem}\label{main}
\rm Let $(\C,\E,\s)$  be an $n$-exangulated category with split idempotents and let $\X$ be a cluster-tilting subcategory of $\C$.
Assume that $\C$ has enough projectives and enough injectives. If $\C$ has Auslander-Reiten $n$-exangles, then $\C/\X$ has Auslander-Reiten $n$-exact sequences.
\end{theorem}

{\bf In order to prove Theorem \ref{main}, we need some preparations as follows.}

\begin{theorem}\label{main11}{\rm\cite[Theorem 3.7]{LZ1}}
\rm Let $\C$ be an $n$-exangulated category with split idempotents and let $\X$ be a cluster-tilting subcategory of $\C$.
 If $\C$ has enough projectives and enough injectives,
then $\C/\X$ is an $n$-abelian category.
\end{theorem}
\begin{remark}
By \cite[Remark 3.9]{LZ1}, we know that the assumption that $\C$ has enough projectives and enough injectives is
necessary in Theorem \ref{main11}.
\end{remark}

\begin{construction}\label{con} Assume that $\X$ is a strongly covariantly finite subcategory of an $n$-exangulated category $\C$.

\textbf{Step 1:} For any object $A\in\C$, take a distinguished $n$-exangle $$A\xrightarrow{f}X_1\xrightarrow{}X_{2}\xrightarrow{}\cdots
\xrightarrow{}X_n\xrightarrow{}X_{n+1}\overset{\delta}{\dashrightarrow},$$ where $f$ is a left $\X$-approximation of $A$ and $X_1,X_2,\cdots,X_{n-1},
X_n\in\X$.

Define $\mathbb{G}(A)=\mathbb{G} A$ to be the image of $X_{n+1}$ in $\C/\X$.

Since $f$ is a left $\X$-approximation of $A$ and $X^{\prime}_1\in\X$, we can complete the following commutative diagram
$$\xymatrix{
A\ar[r]^{f}\ar@{}[dr] \ar[d]^{a} &X_1 \ar[r] \ar@{}[dr]\ar@{-->}[d]^{a_1} &X_2 \ar[r] \ar@{}[dr]\ar@{-->}[d]^{a_2}&\cdot\cdot\cdot \ar[r]\ar@{}[dr] &X_n \ar[r] \ar@{}[dr]\ar@{-->}[d]^{a_n}&X_{n+1} \ar@{}[dr]\ar@{-->}[d]^{b} \ar@{-->}[r]^-{\delta} &\\
{A^{\prime}}\ar[r] &{X^{\prime}_1}\ar[r]&{X^{\prime}_2} \ar[r] &\cdot\cdot\cdot \ar[r] &{X^{\prime} _n}\ar[r]  &{X^{\prime}_{n+1}} \ar@{-->}[r]^-{\delta^\prime} &}
$$
for any morphism $a\in\C(A,A')$.

For any morphism $\bar{a}\in(\C/\X)(A,A')$, define $\mathbb{G}\bar{a}$ to be the image $\bar{b}$ of $b$ in $\C/\X$.

It is not difficult to see that the endofunctor $\mathbb{G}:\C/\X \rightarrow\C/\X$ is well defined.

\textbf{Step 2:} Let $$A_0\xrightarrow{\alpha_0}A_1\xrightarrow{\alpha_1}A_2\xrightarrow{}\cdots\xrightarrow{}A_n\xrightarrow{\alpha_n}A_{n+1}\overset{\eta}{\dashrightarrow}$$ be a distinguished $n$-exangle in $\C$, where $\alpha_0$ is $\X$-monomorphism. Then there exists the following commutative diagram
$$\xymatrix{
A_0\ar[r]^{\alpha_0}\ar@{}[dr] \ar@{=}[d] &A_1 \ar[r]^{\alpha_1} \ar@{}[dr]\ar@{-->}[d]^{g_1} &A_2 \ar[r]^{\alpha_2} \ar@{}[dr]\ar@{-->}[d]^{g_2}&\cdot\cdot\cdot \ar[r]\ar@{}[dr] &A_n \ar[r]^{\alpha_n} \ar@{}[dr]\ar@{-->}[d]^{g_{n}}&A_{n+1} \ar@{}[dr]\ar@{-->}[d]^{{\alpha_{n+1}}} \ar@{-->}[r]^-{\eta} &\\
{A_0}\ar[r]^{f_0} &{X_1}\ar[r]^{f_1}&{X_2} \ar[r]^{f_2} &\cdot\cdot\cdot \ar[r] &{X_n}\ar[r]^{f_n}&\mathbb{G} A_0 \ar@{-->}[r]^-{\zeta} &}
$$
of distinguished $n$-exangles.
Then we have a complex $$A_0\xrightarrow{\overline{\alpha_0}} A_1\xrightarrow{\overline{\alpha_1}}A_2\xrightarrow{\overline{\alpha_2}}\cdots\xrightarrow{\overline{\alpha_n}}A_{n+1}\xrightarrow{\overline{\alpha_{n+1}}}\mathbb{G}A_0.$$ We define right $(n+2)$-angles in $\C/\X$ as the complexes which are isomorphic to complexes obtained in this way.

\end{construction}

\begin{lemma}\label{rig}{\rm\cite[Lemma 3.2]{LZ1}}
\rm {Let $\C$ be an $n$-exangulated category and $\X$ be a strongly covariantly finite subcategory of $\C$. Then the quotient category $\C/\X$ carries an induced right $(n+2)$-angulated structure with the endofunctor $\mathbb{G}$ and right $(n+2)$-angles defined in Construction \ref{con}.}
\end{lemma}
\begin{lemma}\label{seque}{\rm\cite[Lemma 2.7]{ZZ3}}
Let $(\C,\Sigma,\Theta)$  be a right $(n+2)$-angulated category and the sequence $$A_0\xrightarrow{{\alpha_0}} A_1\xrightarrow{{\alpha_1}}A_2\xrightarrow{{\alpha_2}}\cdots\xrightarrow{{\alpha_n}}A_{n+1}\xrightarrow{{\alpha_{n+1}}}\Sigma A_0$$
be a right $(n+2)$-angle. Then $\alpha_i$ is a weak cokernel of $\alpha_{i-1}$ for $i=1,2,\cdots,n+1$.
\end{lemma}

{\bf Now we give the proof of Theorem \ref{main}.}
\proof
By Theorem \ref{main11}, we obtain that $\C/\X$ is an $n$-abelian category.

For any non-projective indecomposable object $C\in\C/\X$, we have $C'\oplus X\in\C$ where $ X\in\X$, and then $C'$ is indecomposable in $\C$ and $C\cong C'$ in $\C/\X$. Since  $\C$ has Auslander-Reiten $n$-exangles, then there exists an Auslander-Reiten $n$-exangle ending at $C'$
\begin{equation}\label{t1}
\begin{array}{l}
A'_0\xrightarrow{\alpha_0}A'_1\xrightarrow{\alpha_1}A'_2\xrightarrow{\alpha_2}\cdots\xrightarrow{\alpha_{n-2}}A'_{n-1}
\xrightarrow{\alpha_{n-1}}A'_n\xrightarrow{\alpha_n}C'\overset{\delta}{\dashrightarrow}.
\end{array}
\end{equation}
Since $\X$ is a cluster-tilting subcategory, there exists a distinguished $n$-exangle
$$X_0\xrightarrow{\beta_0}X_1\xrightarrow{\beta_1}\cdots\xrightarrow{\beta_{n-2}}X_{n-1}\xrightarrow{\beta_{n-1}}X\xrightarrow{\beta}C'\dashrightarrow$$
where $\beta$ is a right $\X$-approximation of $C'$ and $X_0, X_1,\cdots,X\in\X$. We claim that  $\beta$ is not a split epimorphism, otherwise $C'\in \X$ and then $C\simeq C'=0$ in $\C/\X$, this is a contradiction. Since $\alpha_n$ is right almost split, there exists a morphism $\alpha\colon X\to A'_n$ such that $\alpha_n\alpha=\beta$, as the following commutative diagram
$$\xymatrix{& X\ar@{-->}[dl]_{\alpha}\ar[dr]^{\beta} & \\
A'_n\ar[rr]^{\alpha_n} & &C' }$$
By \cite[Proposition 2.14]{HJ}, we know that
\begin{equation}\label{t2}
\begin{array}{l}
0\xrightarrow{}0\xrightarrow{}0\xrightarrow{}\cdots\xrightarrow{}X
\xrightarrow{1}X\xrightarrow{}0\overset{}{\dashrightarrow}
\end{array}
\end{equation}
is a distinguished $n$-exangle. Hence we obtain the following distinguished $n$-exangle  by \cite[Proposition 3.3]{HLN}
$$A'_0\xrightarrow{\alpha_0
}A'_1\xrightarrow{\alpha_1}A'_2\xrightarrow{\alpha_2}\cdots
\xrightarrow{\alpha_{n-3}}A'_{n-2}\xrightarrow{\left(
             \begin{smallmatrix}
              \alpha_{n-2} \\
               0 \\
             \end{smallmatrix}
           \right)
} A'_{n-1}\oplus X\xrightarrow{\left(
             \begin{smallmatrix}
              \alpha_{n-1} & -\alpha\\
               0 & 1 \\
             \end{smallmatrix}
           \right)} A'_{n}\oplus X\xrightarrow{\left(
             \begin{smallmatrix}
              \alpha_{n},  & 0 \\
             \end{smallmatrix}
           \right)} C' \overset{}{\dashrightarrow}.$$
Note that the following commutative diagram
$$\xymatrixcolsep{4.5pc}\xymatrixrowsep{2.7pc}\xymatrix{
A'_0\ar[r]^{\alpha_0
}\ar@{=}[d]& A'_1\ar[r]^{\alpha_1}\ar@{=}[d]& A'_2\ar[r]^{\alpha_2}\ar@{=}[d]& \cdots\ar[r]^{\alpha_{n-2}}&A'_{n-2}\ar@{=}[d]\\
A'_0\ar[r]^{\alpha_0
}& A'_1\ar[r]^{\alpha_1}& A'_2\ar[r]^{\alpha_2}& \cdots\ar[r]^{\alpha_{n-2}}&A'_{n-2}\\
}$$
\begin{equation}\label{t3}
\begin{array}{l}
\xymatrix@C=2cm{
\ar[r]^{\left(
             \begin{smallmatrix}
              \alpha_{n-2} \\
               0 \\
             \end{smallmatrix}
           \right)\qquad} & A'_{n-1}\oplus X\ar[r]^{\left(
             \begin{smallmatrix}
              \alpha_{n-1} & -\alpha\\
               0 & 1 \\
             \end{smallmatrix}
           \right) }\ar@{=}[d]&  A'_{n}\oplus X\ar[r]^{\left(
             \begin{smallmatrix}
              \alpha_{n}, & \beta \\
             \end{smallmatrix}
           \right)}\ar[d]_{\cong}^{\left(
             \begin{smallmatrix}
            1 & \alpha\\
               0 & 1 \\
             \end{smallmatrix}
           \right)\ \ \ \ }&  C'\overset{}{\dashrightarrow}\hspace{-7mm}\ar@{=}[d]\\
\ar[r]^{\left(
             \begin{smallmatrix}
              \alpha_{n-2} \\
               0 \\
             \end{smallmatrix}
           \right)\qquad}& A'_{n-1}\oplus X\ar[r]^{\left(
             \begin{smallmatrix}
              \alpha_{n-1} & 0\\
               0 & 1 \\
             \end{smallmatrix}
           \right)\ }& A'_{n}\oplus X\ar[r]^{\left(
             \begin{smallmatrix}
               \alpha_{n}, & 0 \\
             \end{smallmatrix}
           \right)}&C'\overset{}{\dashrightarrow}\hspace{-7mm}}
\end{array}
\end{equation}
which shows that the first row of (\ref{t3}) is a distinguished $n$-exangle. It is easy to see that $(\alpha_{n},\beta)$ is $\X$-epic since $\beta$ is $\X$-epic. Moreover, we have the commutative diagram
$$\xymatrixcolsep{4.5pc}\xymatrixrowsep{2.7pc}\xymatrix{
X_0\ar[r]^{\beta_0
}\ar@{-->}[d]^{f_0}& X_1\ar[r]^{\beta_1}\ar@{-->}[d]^{f_1}& X_2\ar[r]^{\beta_2}\ar@{-->}[d]^{f_2}& \cdots\ar[r]^{\beta_{n-3}}&  X_{n-2}\ar@{-->}[d]^{f_{n-2}}\\
A'_0\ar[r]^{\alpha_0
}& A'_1\ar[r]^{\alpha_1}& A'_2\ar[r]^{\alpha_2}& \cdots\ar[r]^{\alpha_{n-3}}&  A'_{n-2}\\
}$$
\begin{equation}\label{t4}
\begin{array}{l}
\xymatrix@C=2cm{
\ar[r]^{\beta_{n-2}}&  X_{n-1}\ar[r]^{\beta_{n-1} }\ar@{-->}[d]^{f_{n-1}}&  X\ar[r]^{\beta}\ar[d]^{\left(
             \begin{smallmatrix}
           0\\
               1 \\
             \end{smallmatrix}
           \right)\ \ \ \ }&  C'\overset{}{\dashrightarrow}\hspace{-7mm}\ar@{=}[d]\\
\ar[r]^{\left(
             \begin{smallmatrix}
              \alpha_{n-2} \\
               0 \\
             \end{smallmatrix}
           \right)}& A'_{n-1}\oplus X\ar[r]^{\left(
             \begin{smallmatrix}
              \alpha_{n-1} & -\alpha\\
               0 & 1 \\
             \end{smallmatrix}
           \right)\ }& A'_{n}\oplus X\ar[r]^{\left(
             \begin{smallmatrix}
               \alpha_{n} & \beta \\
             \end{smallmatrix}
           \right)}& C'\overset{}{\dashrightarrow}\hspace{-7mm}}
\end{array}
\end{equation}
of distinguished $n$-exangles by the dual of Proposition 3.6 in \cite{HLN}. So we obtain a left $(n+2)$-angle by the dual of Lemma \ref{rig}

\begin{equation}\label{t5}
\begin{array}{l}0\xrightarrow{}A_0\xrightarrow{\overline{\alpha_0}} A_1\xrightarrow{\overline{\alpha_1}}A_2\xrightarrow{\overline{\alpha_2}}\cdots\xrightarrow{\overline{\alpha_{n-2}}}A_{n-1}\xrightarrow{\overline{\alpha_{n-1}}}A_{n}\xrightarrow{\overline{\alpha_{n}}} C.
\end{array}
\end{equation}

From the dual of Lemma \ref{seque}, we know that (\ref{t5}) is a left $n$-exact sequence.

Since $\X$ is a cluster-tilting subcategory, there exists a distinguished $n$-exangle
$$A'_0\xrightarrow{\gamma}Y\xrightarrow{\gamma_1}Y_2\xrightarrow{\gamma_2}\cdots\xrightarrow{\gamma_{n-2}}Y_{n-1}\xrightarrow{\gamma_{n-1}}Y_{n}\xrightarrow{\gamma_{n}}Y_{n+1}\dashrightarrow$$
where $\gamma$ be a left $\X$-approximation of $A'_0$ and $Y, Y_2,\cdots,Y_{n+1}\in\X$. Note that $\gamma$ is not a split monomorphism, since $\alpha_0$ is left almost split, there exists a morphism $l\colon A'_1\to Y$ such that $l\alpha_0=\gamma$, as the following commutative diagram
$$\xymatrix@R=0.5cm{A'_0\ar[dr]^{\gamma} \ar[rr]^{\alpha_0} & & A'_1\ar@{-->}[dl]_{l}\\
& Y &}$$
By \cite[Proposition 2.14]{HJ}, we know that
\begin{equation}\label{t6}
\begin{array}{l}
0\xrightarrow{}Y\xrightarrow{1}Y\xrightarrow{}0\xrightarrow{}\cdots
\xrightarrow{}0\xrightarrow{}0\overset{}{\dashrightarrow}
\end{array}
\end{equation}
is a distinguished $n$-exangle. By \cite[Proposition 3.3]{HLN},  we obtain the following distinguished $n$-exangle
$$\xymatrixcolsep{4.5pc} \xymatrix{
A'_0\ar[r]^{\left(
             \begin{smallmatrix}
              \alpha_{0} \\
               0 \\
             \end{smallmatrix}
           \right)
}& A'_1\oplus Y\ar[r]^{\left(
             \begin{smallmatrix}
              \alpha_{1} & 0\\
               0 & 1 \\
             \end{smallmatrix}
           \right)\ \ }& A'_2\oplus Y\ar[r]^{\ \ \ \left(
             \begin{smallmatrix}
              \alpha_{2},  & 0\\
             \end{smallmatrix}
           \right)}&  A'_3\ar[r]^{\alpha_3}&\cdots\\
}$$
$$\begin{gathered}\xymatrixcolsep{5pc}\xymatrix{
\cdots\ar[r]^{\alpha_{n-3}}& A'_{n-2} \ar[r]^{\left(
             \begin{smallmatrix}
              \alpha_{n-2} \\
               0 \\
             \end{smallmatrix}
           \right)
}& A'_{n-1}\oplus X\ar[r]^{\left(
             \begin{smallmatrix}
              \alpha_{n-1} & -\alpha\\
               0 & 1 \\
             \end{smallmatrix}
           \right)\ \ \ \ }& A'_{n}\oplus X\ar[r]^{\ \ \ \left(
             \begin{smallmatrix}
              \alpha_{n},  & \beta\\
             \end{smallmatrix}
           \right)}& C' \overset{}{\dashrightarrow},
}\end{gathered}$$
which is the sum to both the first row of (\ref{t3}) and (\ref{t6}).
Note that the following commutative diagram
$$\xymatrixcolsep{4.5pc}\xymatrixrowsep{2.7pc}\xymatrix{
A'_0\ar[r]^{\left(
             \begin{smallmatrix}
              \alpha_{0} \\
              \gamma \\
             \end{smallmatrix}
           \right)
}\ar@{=}[d]& A'_1\oplus Y\ar[r]^{\left(
             \begin{smallmatrix}
             \alpha_1 & 0\\
               -l & 1 \\
             \end{smallmatrix}
           \right)\ \\\\ }\ar[d]_{\cong}^{\left(
             \begin{smallmatrix}
            1 &0\\
               -l & 1 \\
             \end{smallmatrix}
           \right)\ \ \ \ }& A'_2\oplus Y\ar[r]^{\left(
             \begin{smallmatrix}
               \alpha_{2}, & 0 \\
             \end{smallmatrix}
           \right)}\ar@{=}[d]& A'_3\ar[r]^{\alpha_3
}\ar@{=}[d]& \cdots\ar[r]^{\alpha_{n-3}}& A'_{n-2}\ar@{=}[d]\\
A'_0\ar[r]^{\left(
             \begin{smallmatrix}
              \alpha_{0} \\
               0 \\
             \end{smallmatrix}
           \right)
}& A'_1\oplus Y\ar[r]^{\left(
             \begin{smallmatrix}
             \alpha_1 & 0\\
               0& 1 \\
             \end{smallmatrix}
           \right)\ \\\\ }& A'_2\oplus Y\ar[r]^{\left(
             \begin{smallmatrix}
               \alpha_{2}, & 0 \\
             \end{smallmatrix}
           \right)}& A'_3\ar[r]^{\alpha_3}& \cdots\ar[r]^{\alpha_{n-3}}& A'_{n-2}\\
}$$
\begin{equation}\label{t7}
\begin{array}{l}
\xymatrix@C=2cm{
\ar[r]^{\left(
             \begin{smallmatrix}
              \alpha_{n-2} \\
               0 \\
             \end{smallmatrix}
           \right)} & A'_{n-1}\oplus X\ar[r]^{\left(
             \begin{smallmatrix}
              \alpha_{n-1} & -\alpha\\
               0 & 1 \\
             \end{smallmatrix}
           \right)\ }\ar@{=}[d]&  A'_{n}\oplus X\ar[r]^{\left(
             \begin{smallmatrix}
              \alpha_{n}, & \beta \\
             \end{smallmatrix}
           \right)}\ar@{=}[d]&  C'\overset{}{\dashrightarrow}\hspace{-7mm}\ar@{=}[d] \\
\ar[r]^{\left(
             \begin{smallmatrix}
              \alpha_{n-2} \\
               0 \\
             \end{smallmatrix}
           \right)}& A'_{n-1}\oplus X\ar[r]^{\left(
             \begin{smallmatrix}
              \alpha_{n-1} & -\alpha\\
               0 & 1 \\
             \end{smallmatrix}
           \right)\ }& A'_{n}\oplus X\ar[r]^{\left(
             \begin{smallmatrix}
               \alpha_{n}, & \beta \\
             \end{smallmatrix}
           \right)}& C'\overset{}{\dashrightarrow}\hspace{-7mm}}
\end{array}
\end{equation}
which shows that the first row of (\ref{t7}) is a distinguished $n$-exangle. It is easy to see that $\left(
             \begin{smallmatrix}
              \alpha_{0} \\
              \gamma \\
             \end{smallmatrix}
           \right)$ is $\X$-monic since $\gamma$ is $\X$-monic. Moreover, we have the commutative diagram
$$\xymatrixcolsep{4.5pc}\xymatrixrowsep{2.7pc}\xymatrix{
A'_0\ar[r]^{\left(
             \begin{smallmatrix}
              \alpha_{0} \\
              \gamma \\
             \end{smallmatrix}
           \right)
}\ar@{=}[d]& A'_1\oplus Y\ar[r]^{\left(
             \begin{smallmatrix}
             \alpha_1 & 0\\
               -l & 1 \\
             \end{smallmatrix}
           \right)\ \\\\ }\ar[d]^{\left(
             \begin{smallmatrix}
               0, & 1\\
             \end{smallmatrix}
           \right) }& A'_2\oplus Y\ar[r]^{\left(
             \begin{smallmatrix}
               \alpha_{2}, & 0 \\
             \end{smallmatrix}
           \right)}\ar@{-->}[d]^{h_2}& A'_3\ar[r]^{\alpha_3
}\ar@{-->}[d]^{h_3}& \cdots\ar[r]^{\alpha_{n-3}}& A'_{n-2}\ar@{-->}[d]^{h_{n-2}}\\
A'_0\ar[r]^{\gamma
}& Y\ar[r]^{\gamma_1}& Y_2\ar[r]^{\gamma_2}& Y_3\ar[r]^{\gamma_3}& \cdots\ar[r]^{\gamma_{n-3}}& Y_{n-2}\\
}$$
\begin{equation}\label{t8}
\begin{array}{l}
\xymatrix@C=2cm{
 \ar[r]^{\left(
             \begin{smallmatrix}
              \alpha_{n-2} \\
               0 \\
             \end{smallmatrix}
           \right)} & A'_{n-1}\oplus X\ar[r]^{\left(
             \begin{smallmatrix}
              \alpha_{n-1} & -\alpha\\
               0 & 1 \\
             \end{smallmatrix}
           \right)\ }\ar@{-->}[d]^{h_{n-1}}&  A'_{n}\oplus X\ar[r]^{\left(
             \begin{smallmatrix}
              \alpha_{n}, & \beta \\
             \end{smallmatrix}
           \right)}\ar@{-->}[d]^{h_{n}}&   C'\overset{}{\dashrightarrow}\hspace{-7mm}\ar@{-->}[d]^{h_{n+1}} \\
  \ar[r]^{\gamma_{n-2}}& Y_{n-1}\ar[r]^{\gamma_{n-1}}&Y_{n}\ar[r]^{\gamma_{n}}& Y_{n+1}\overset{}{\dashrightarrow}\hspace{-7mm}}
\end{array}
\end{equation}
of distinguished $n$-exangles by Proposition 3.6 in \cite{HLN}. So we obtain a right $(n+2)$-angle by Lemma \ref{rig}
\begin{equation}\label{t9}
\begin{array}{l}
A_0\xrightarrow{\overline{\alpha_0}} A_1\xrightarrow{\overline{\alpha_1}}A_2\xrightarrow{\overline{\alpha_2}}\cdots\xrightarrow{\overline{\alpha_{n-2}}}A_{n-1}\xrightarrow{\overline{\alpha_{n-1}}}A_{n}\xrightarrow{\overline{\alpha_{n}}} C\xrightarrow{}0.
\end{array}
\end{equation}
By Lemma \ref{seque}, we get that (\ref{t5})  is a right $n$-exact sequence. Thus \begin{equation}\label{t10}
\begin{array}{l}0\xrightarrow{} A_0\xrightarrow{\overline{\alpha_0}} A_1\xrightarrow{\overline{\alpha_1}}A_2\xrightarrow{\overline{\alpha_2}}\cdots\xrightarrow{\overline{\alpha_{n-2}}}A_{n-1}\xrightarrow{\overline{\alpha_{n-1}}}A_{n}\xrightarrow{\overline{\alpha_{n}}} C\xrightarrow{}0.\end{array}
\end{equation}
is a $n$-exact sequence.

Finally, we prove that the above $n$-exact sequence (\ref{t10}) is an Auslander-Reiten $n$-exact sequence. Since $\alpha_n\colon A'_n\to C'$ is right almost split, $\alpha_n$ is not a split epimorphism, that is to say $C'$ is not a direct summand of $ A'_n$. Suppose that ${\overline{\alpha_n}}\colon A_n\to C$ is a split epimorphism, then there exists a morphism ${\overline{d_n}}\colon C\to A_n$, such that ${\overline{\alpha_n}}\circ{\overline{d_n}}= 1_C$. So there is an object $X\in\X$, such that $1_{C'}-\alpha_nd_n=ts$, it reflected in the following commutative diagram
$$\xymatrix{C' \ar[rr]^{1_{C'}-\alpha_nd_n}\ar[dr]^{s} & & C'\\
& X \ar[ur]^{t}}$$
Then $(\alpha_{n},t)\left(
             \begin{smallmatrix}
              d_{n} \\
              s \\
             \end{smallmatrix}
           \right)=\alpha_nd_n+ts=1_{C'}$, that is, $(\alpha_{n},t)\colon A'_n\oplus X\to C'$ is a split epimorphism. Thus $C'$ is a direct summand of $ A'_n\oplus X$. Since $\alpha_n\colon A'_n\to C'$ is right almost split, then $C'$ is indecomposable in $\C$, in particular, $C'$ is not a direct summand of $ X$. Thus $C'$ is a direct summand of $ A'_n$, this is a contradiction. This shows that ${\overline{\alpha_n}}$ is not a split epimorphism. For every object $D\in\C/\X$, every morphism ${\overline{p_n}}\colon D\to C$, which is not a split epimorphism, then ${{p_n}}\colon D'\to C'$ is not a split epimorphism, there exists a morphism ${{q_n}}\colon D'\to  A'_n$ such that $\alpha_nq_n=p_n$, it reflected in the following commutative diagram
           $$\xymatrix{&  D'\ar@{-->}[dl]_{q_n}\ar[dr]^{p_n} & \\
A'_n\ar[rr]^{\alpha_n} & &C' }$$
Then we have ${\overline{\alpha_n}}\circ{\overline {q_n}}={\overline{ p_n}}$, and hence ${\overline{\alpha_n}}$ is right almost split. Similarly, we can prove that ${\overline{\alpha_0}}$ is left almost split.

Since (\ref{t1}) is an Auslander-Reiten $n$-exangle, then $\alpha_i\colon A'_i\to A'_{i+1}\in {\rm rad}_{\C}$, where $i=1,2,\cdots,n-1$. Thus for any morphism ${\overline {h_i}}\colon A_{i+1}\to A_{i}$, we have  ${h_i}\colon A'_{i+1}\to A'_{i}\in \C$ such that $1_{A'_{i}}-h_{i}\alpha_{i}$ is invertible. And hence there exists a morphism $\varepsilon_{i}\colon A'_{i}\to A'_{i}$, such that
$\varepsilon_{i}(1_{A'_{i}}-h_{i}\alpha_{i})=1_{A'_{i}}$ and $(1_{A'_{i}}-h_{i}\alpha_{i})\varepsilon_{i}=1_{A'_{i}}$ hold. Then we have $\overline {\varepsilon_{i}}(\overline {1}_{A'_{i}}-\overline {h_{i}}\overline {\alpha_{i}})=\overline {1}_{A'_{i}}$ and $(\overline {1}_{A'_{i}}-\overline {h_{i}}\overline {\alpha_{i}})\overline {\varepsilon_{i}}=\overline {1}_{A'_{i}}$. Therefore, $\overline {\alpha_i}\in {\rm rad}_{\C/\X}$, where $i=1,2,\cdots,n-1$. This shows that the $n$-exact sequence (\ref{t10}) is an Auslander-Reiten $n$-exact sequence.

Dually, for any non-injective indecomposable object $C\in\C/\X$, we can get an Auslander-Reiten $n$-exact sequence starting at $C$:
$$0\xrightarrow{} C\xrightarrow{\overline{\varrho_0}} B_1\xrightarrow{\overline{\varrho_1}}B_2\xrightarrow{\overline{\varrho_2}}\cdots\xrightarrow{\overline{\varrho_{n-2}}}B_{n-1}\xrightarrow{\overline{\varrho_{n-1}}}B_{n}\xrightarrow{\overline{\varrho_{n}}} B_{n+1}\xrightarrow{}0.$$
Therefore, $\C/\X$ has Auslander-Reiten $n$-exact sequences.\qed

By applying Theorem \ref{main} to $(n+2)$-angulated categories, and using the fact that any
$(n+2)$-angulated category can be viewed as an $n$-exangulated category with enough
projectives and injectives, we get the following result.
\begin{corollary}{\rm \cite[Theorem 3.4]{LH}}
\rm \rm Let $\C$ be an $(n+2)$-angulated category with split idempotents and $\X$ a cluster-tilting subcategory of $\C$.
If $\C$ has Auslander-Reiten $(n+2)$-angles, then $\C/\X$ has Auslander-Reiten $n$-exact sequences.
\end{corollary}

\begin{corollary}
\rm Let $(\C,\E,\s)$  be an $k$-linear Hom-finite Krull-Schmidt $n$-exangulated category with split idempotents and $\X$ a cluster-tilting subcategory of $\C$.
Assume that $\C$ has enough projectives and enough injectives. If $\C$ is locally finite, then $\C/\X$ has Auslander-Reiten $n$-exact sequences.
\end{corollary}

\proof
This follows that \cite[Theorem 3.12]{HHZZ} and Theorem \ref{main}.
\qed

\begin{definition}\cite[Definition 3.2]{LZ} and \cite[Definition 3.17]{ZW}
 Let $(\C,\E,\s)$  be an $n$-exangulated category. $\C$ is said to be Frobenius if $\C$ has enough projectives and  enough injectives and if moreover the projectives coincide with the injectives.
 \end{definition}

Let $(\C,\E,\s)$ be a Frobenius $n$-exangulated category. We denote its stable category by $\overline{\C}$. Our second main result establishes a connection between Auslander-Reiten $n$-exangles in Frobenius $n$-exangulated categories and Auslander-Reiten $(n+2)$-angles in $(n+2)$-angulated categories.

\begin{theorem}\label{main1}
\rm  Let $(\C,\E,\s)$  be a Frobenius $n$-exangulated category. If $\C$ has Auslander-Reiten $n$-exangles, then the stable category $\overline{\C}$ of $\C$ has Auslander-Reiten $(n+2)$-angles.
\end{theorem}
{\bf In order to prove Theorem \ref{main1}, we need some preparations as follows.}

\begin{construction}\label{conj} Assume that $(\C,\E,\s)$  is a Frobenius $n$-exangulated category
and we denote the full subcategory of injective objects in $\C$ by $\I$.

\textbf{Step 1:} Any object $X\in\C$ admits a distinguished $n$-exangle $$X\xrightarrow{i_0}I_1\xrightarrow{i_1}I_{2}\xrightarrow{i_2}\cdots
\xrightarrow{i_{n-1}}I_n\xrightarrow{i_n}Y\overset{\delta_X}{\dashrightarrow},$$ where $i_0$ is a left $\I$-approximation of $X$, $i_n$ is a right $\I$-approximation of $Y$ and $I_i\in\I$ for $i\in\{1,2,\cdots,n\}$.

Define $S(X)=S X$ to be the image of $Y$ in $\overline\C$. Then the functor $S: \overline\C\rightarrow\overline\C$ is an automprhism by \cite[Proposition 3.7]{LZ} or
\cite[Proposition 3.17]{ZW}.

\textbf{Step 2:} For $\overline \C$, let $\overline \E\colon\overline{\C}^{\op} \times \overline \C\to \Ab$ and $\overline {\mathfrak{s}}$ be the bifunctor given by
\begin{itemize}
\item $\overline \E(C,A)=\E(C,A)$, $\forall A,C\in \C$.
\item $\overline \E(\overline c, \overline a)=\E(c,a)$, $\forall a\in \E(A,A'),c\in \E(C,C')$.
\item For any $\overline \E$-extension $\delta\in \overline \E(C,A)=\E(C,A)$, define
$$\overline {\mathfrak{s}}(\delta)=\overline {\mathfrak{s}(\delta)}=[\xymatrix{A \ar[r]^{\overline {d_0}} &X_1 \ar[r]^{\overline {d_1}} \ar[r] &\cdot\cdot\cdot \ar[r]^{\overline {d_{n-1}}} &X_{n} \ar[r]^{\overline {d_{n}}} &C \ar@{-->}[r]^-{\delta} &].}$$
\end{itemize}
By \cite[Lemma 3.12]{LZ}, we know that $(\C,\E,\s)$ is an $n$-exangulated category. Moreover, we have a one-to-one correspondence between $\E(C,A)$ and $\overline {\C}(C,SA)$ by \cite[Lemma 3.11]{LZ}. So for any distinguished $n$-exangle
$$A\xrightarrow{\overline {d_0}}X_1\xrightarrow{\overline {d_1}}X_{2}\xrightarrow{\overline {d_2}}\cdots
\xrightarrow{\overline {d_{n-1}}}X_n\xrightarrow{\overline {d_n}}C\overset{\delta}{\dashrightarrow},$$
it induces an $(n+2)$-$S$-sequence
$$A\xrightarrow{\overline {d_0}}X_1\xrightarrow{\overline {d_1}}X_{2}\xrightarrow{\overline {d_2}}\cdots
\xrightarrow{\overline {d_{n-1}}}X_n\xrightarrow{\overline {d_n}}C\xrightarrow{\overline {d_{n+1}}}SA.$$
Let $\Theta$ be the class of $(n+2)$-$S$-sequences. Then we have the following.
\end{construction}

\begin{lemma}\label{rp}{\rm\cite[Theorem 3.13]{LZ}}
\rm {Let  $(\C,\E,\s)$ be a Frobenius $n$-exangulated category. Then $(\overline{\C},S,\Theta)$ is an $(n+2)$-angulated category.}
\end{lemma}

{\bf Now we give the proof of Theorem \ref{main1}.}
\proof
By Lemma \ref{rp}, we know that $\overline{\C}$ is an $(n+2)$-angulated category.

For any indecomposable object $C\in\overline{\C}$, we have $C'\oplus I\in\C$ where $ I\in\I$, and then $C'$ is indecomposable in $\C$ and $C\cong C'$ in $\overline{\C}$. Since  $\C$ has Auslander-Reiten $n$-exangles, then there exists an Auslander-Reiten $n$-exangle ending at $C'$
\begin{equation}\label{t1}
\begin{array}{l}
A'_0\xrightarrow{\alpha_0}A'_1\xrightarrow{\alpha_1}A'_2\xrightarrow{\alpha_2}\cdots\xrightarrow{\alpha_{n-2}}A'_{n-1}
\xrightarrow{\alpha_{n-1}}A'_n\xrightarrow{\alpha_n}C'\overset{\delta}{\dashrightarrow}.
\end{array}
\end{equation}
By Construction \ref{conj}, we obtain an $(n+2)$-angle in $\overline{\C}$
\begin{equation}\label{tm}
\begin{array}{l}
A_0\xrightarrow{\overline{\alpha_0}}A_1\xrightarrow{\overline{\alpha_1}}A_2\xrightarrow{\overline{\alpha_2}}\cdots\xrightarrow{\overline{\alpha_{n-2}}}A_{n-1}
\xrightarrow{\overline{\alpha_{n-1}}}A_n\xrightarrow{\overline{\alpha_n}}C\xrightarrow{\overline{\alpha_{n+1}}}SA_0.
\end{array}
\end{equation}
Using similar arguments as in the proof of Theorem \ref{main}, we get that $\overline{\alpha_0}$ is left almost split and $\overline{\alpha_n}$ is right almost split. Moreover, $\overline{\alpha_1},\overline{\alpha_2},\cdots, \overline{\alpha_{n-1}}\in{\rm rad}_{\overline{\C}}$.  This shows that the $(n+2)$-angle \ref{tm} is an Auslander-Reiten $(n+2)$-angle.

Dually, for any indecomposable object $C\in{\overline\C}$, we can get an Auslander-Reiten$(n+2)$-angle starting at $C$:
$$ C\xrightarrow{\overline{\varrho_0}} B_1\xrightarrow{\overline{\varrho_1}}B_2\xrightarrow{\overline{\varrho_2}}\cdots\xrightarrow{\overline{\varrho_{n-2}}}B_{n-1}\xrightarrow{\overline{\varrho_{n-1}}}B_{n}\xrightarrow{\overline{\varrho_{n}}} B_{n+1}\xrightarrow{\overline{\varrho_{n+1}}}SC.$$
Therefore, $\overline{\C}$ has Auslander-Reiten $(n+2)$-angles.
\qed

\vspace{2mm}

In Theorem \ref{main1}, when $n=1$, we have the following.
\begin{corollary}
\rm  Let $(\C,\E,\s)$  be a Frobenius extriangulated category. If $\C$ has Auslander-Reiten $\E$-triangles, then the stable category $\overline{\C}$ of $\C$ has Auslander-Reiten triangles.
\end{corollary}

\vspace{0.5cm}

%

\textbf{Jian He}\\
Department of Applied Mathematics, Lanzhou University of Technology, 730050 Lanzhou, Gansu, P. R. China\\
E-mail: \textsf{jianhe30@163.com}\\[0.3cm]
\textbf{Hangyu Yin}\\
College of Mathematics, Hunan Institute of Science and Technology, 414006 Yueyang, Hunan, P. R. China.\\
E-mail: hangyuyin@163.com\\[0.3cm]
\textbf{Panyue Zhou}\\
School of Mathematics and Statistics, Changsha University of Science and Technology, 410114 Changsha, Hunan,  P. R. China\\
E-mail: \textsf{panyuezhou@163.com}

\end{document}